\documentclass[12pt,letterpaper,reqno]{amsart}
\usepackage{mathrsfs}
\usepackage{times}
\usepackage[T1]{fontenc}
\usepackage{amsmath,amsfonts,amsthm,amssymb,amscd}
\usepackage[dvips]{graphics}
\usepackage{epsfig}

\newcommand{\ncr}[2]{{#1 \choose #2}}
\newcommand{\sumin}{\sum_{i=1}^N}
\newcommand{\fon}{\frac{1}{N}}

\newcommand{\glik}{\gl_i^k(A)}
\newcommand{\mtk}{M_{2k}(N)}
\newcommand{\foh}{\frac{1}{2}}

\newtheorem{thm}{Theorem}[section]

\newtheorem{lem}[thm]{Lemma}

\newtheorem{defi}[thm]{Definition}

\newtheorem{rek}[thm]{Remark}

\newcommand\bd{\begin{displaymath}}
\newcommand\ed{\end{displaymath}}
\newcommand\be{\begin{equation}}
\newcommand\ee{\end{equation}}
\newcommand\bea{\begin{eqnarray}}
\newcommand\eea{\end{eqnarray}}
\newcommand\bi{\begin{itemize}}
\newcommand\ei{\end{itemize}}
\newcommand\ben{\begin{enumerate}}
\newcommand\een{\end{enumerate}}
\newcommand\bc{\begin{center}}
\newcommand\ec{\end{center}}
\newcommand\ba{\begin{array}}
\newcommand\ea{\end{array}}

\newcommand{\Inv}[1]{\frac{1}{#1}}
\newcommand{\bs}[2]{b_{|i_{#1}-i_{#2}|}}

\newcommand{\R}{\ensuremath{\Bbb{R}}}

\newcommand{\gep}{\epsilon}   
\newcommand{\gl}{\lambda}

\newcommand{\E}{\ensuremath{\mathbb{E}}}

\newcommand{\N}{\mathbb{N}}

\begin{document}

\title[Eigenvalue Distribution for Real
Symmetric Toeplitz Matrices]{Eigenvalue Spacing Distribution for
the Ensemble of Real Symmetric Toeplitz Matrices}

\email{hammond.162@osu.edu, sjmiller@math.ohio-state.edu}

\author{Christopher Hammond, Steven J. Miller}

\address{Department of Mathematics\\ The Ohio State University\\
Columbus, OH 43210 U.S.A.}

\thanks{This work was done at the American Institute for Mathematics
and The Ohio State University. It is a pleasure to thank them for
their help and support. We are grateful to the participants of
AIM's Summer 2003 REU program, Boris Pittel, and Eitan Sayag for
helpful comments and discussions, and to Wlodzimierz Bryc, Amir
Dembo and Tiefeng Jiang for sharing their preprint. We would
especially like to thank Peter March for showing us how if we
could prove something along the lines of Theorem
\ref{thmprestrong}, then almost sure convergence would follow.}

\begin{abstract}
Consider the ensemble of Real Symmetric Toeplitz Matrices, each
entry iidrv from a fixed probability distribution p of mean 0,
variance 1, and finite higher moments. The limiting spectral
measure (the density of normalized eigenvalues) converges weakly
to a new universal distribution with unbounded support,
independent of p. This distribution's moments are almost those of
the Gaussian's; the deficit may be interpreted in terms of
Diophantine obstructions. With a little more work, we obtain
almost sure convergence. An investigation of spacings between
adjacent normalized eigenvalues looks Poissonian, and not GOE.\\

Classification: 15A52 (primary), 60F99, 62H10 (secondary).\\

Keywords: Random Matrix Theory, Toeplitz Matrices, Distribution of
Eigenvalues, Diophantine Obstructions
\end{abstract}

\maketitle

\section{Introduction}

One of the central problems in Random Matrix Theory is as follows:
consider some ensemble of matrices $A$ with probabilities $p(A)$.
As $N \to \infty$, what can one say about the density of
normalized eigenvalues? For Real Symmetric matrices, where the
entries are iidrv from suitably restricted probability
distributions, the limiting distribution is the semi-circle. Note
this ensemble has $\frac{N(N+1)}{2}$ independent parameters
($a_{ij}, i \le j$). For matrix ensembles with fewer degrees of
freedom, different limiting distributions arise (for example,
McKay \cite{McK} proved $d$-regular graphs are given by Kesten's
Measure). By examining ensembles with fewer than $N^2$ degrees of
freedom, one has the exciting potential of seeing new, universal
distributions. In this paper we investigate Symmetric Toeplitz
matrices.

\begin{defi}

A Toeplitz matrix is a matrix of the form

\be \left(\begin{array}{ccccc}
                        b_{0}  &  b_{1}  & b_{2}  & \cdots & b_{N-1} \\
                        b_{-1} &  b_{0}  & b_{1}  & \cdots & b_{N-2} \\
                        b_{-2} &  b_{-1} & b_{0}  & \cdots & b_{N-3} \\
                        \vdots & \vdots  & \vdots & \ddots & \vdots \\
            b_{1-N} &  b_{2-N} & b_{3-N} & \cdots & b_0 \\
                          \end{array}\right)
\ee
\end{defi}

We investigate symmetric Toeplitz matrices whose entries are
chosen according to some distribution $p$ with mean 0, variance 1,
and finite higher moments. The probability density of a given
matrix is $\prod_{i=0}^{N-1} p(b_i)$.

By looking at $\mbox{Trace}(A^2) = \sum_i \gl_i^2(A)$, we see that
the eigenvalues of $A$ are of order $\sqrt{N}$. As the main
diagonal is constant, all $b_0$ does is shift each eigenvalue.
Therefore, it is sufficient to consider the case where the main
diagonal vanishes.

To each Toeplitz matrix, we may attach a spacing measure by
placing a point mass of size $\frac{1}{N}$ at each normalized
eigenvalue:

\be \mu_{A,N}(x)dx \ = \ \frac{1}{N} \sum_{i=1}^N \delta\left( x -
\frac{\gl_i(A)}{\sqrt{N}} \right)dx. \ee

The $k^{th}$ moment of $\mu_{A,N}(x)$ is

\be M_k(A,N) \ = \ \frac{1}{N^{\frac{k}{2}+1}} \sum_{i=1}^N
\gl_i^k(A). \ee

Let $M_k(N)$ be the average of $M_k(A,N)$ over the ensemble, with
each $A$ weighted by its density. We show that $M_k(N)$ converges
to the moments of a new universal distribution, independent of
$p$. The new distribution looks Gaussian, and numerical
simulations and heuristic sketches at first seemed to support such
a conjecture. A more detailed analysis, however, reveals that
while $M_k(N)$ agrees with the Gaussian moments for odd $k$ and $k
= 0, 2$, the other even moments are less than the Gaussian.

We now sketch the proof. By the Trace Lemma,

\be \sumin \glik \ = \ \mbox{Trace}(A^k) \ = \ \sum_{1 \leq
i_{1},\dotsm, i_{k} \leq N}a_{i_{1},i_{2}}a_{i_{2},i_{3}}\cdots
a_{i_{k},i_{1}}. \ee

As our Toeplitz matrices are constant along diagonals, depending
only on $|i_m - i_n|$, we have

\be\label{eqtr} M_k(N) \ = \ \frac{1}{N^{\frac{k}{2}+1}}  \sum_{1
\leq i_{1},\cdots, i_k \leq
N}\E(b_{|i_{1}-i_{2}|}b_{|i_{2}-i_{3}|}\dotsm b_{|i_{k}-i_{1}|}),
\ee

where by $\E(\cdots)$ we mean averaging over the Toeplitz
ensemble, with each matrix $A$ weighted by its probability of
occurring, and the $b_j$ are iidrv drawn from $p(x)$.

We then show that as $N \to \infty$, the above sums vanish for $k$
odd, and converge independent of $p$ for $k$ even to numbers $M_k$
bounded by the moments of the Gaussian. By showing $\E[ |M_k(A,N)
- M_k(N)|^m]$ is small for $m=2$ ($m=4$), we obtain weak (almost
sure) convergence.

\begin{rek} This problem was first posed by Bai \cite{Bai}, where
he also asked similar questions about Hankel and Markov matrices.
Almost surely the methods of this paper would be applicable to
these cases. Bose and Bryc-Dembo-Jiang have independently observed
that the limiting distribution is not Gaussian. Using a more
probabilistic formulation, \cite{BDJ} have calculated the moments
using uniform variables and interpreting results as volumes of
solids related to Eulerian numbers. We have independently found
the same numbers, but through Diophantine analysis, which allows
us to interpret the deviations from the Gaussian in terms of
Diophantine obstructions, and estimate the rate of convergence.
\end{rek}

\section{Determination of the Moments}

\subsection{$k = 0, 2$ and $k$ odd}

For all $N$, $M_0(A,N) = M_0(N) = 1$. For $k = 2$, we have

\be M_2(N) \ = \ \frac{1}{N^2} \sum_{1 \le i_1,i_2 \le N} \E(
b_{|i_1-i_2|} b_{|i_2 - i_1}) \ = \ \frac{1}{N^2} \sum_{1 \le
i_1,i_2 \le N} \E(b_{|i_1-i_2|}^2). \ee

As we have drawn the $b$s from a variance one distribution, the
expected value above is $1$ if $i_1 \neq i_2$ and $0$ otherwise.
Thus, $M_2(N) = \frac{N^2 - N}{N^2} = 1 - \fon$. Note there are
two degrees of freedom. We can choose $b_{|i_1 - i_2|}$ to be on
any diagonal. Once we have specified the diagonal, we can then
choose $i_1$ freely, which now determines $i_2$.

For $k$ odd, we must have at least one $b_j$ occurring to an odd
power. If one occurs to the first power, as the expected value of
a product of independent variables is the product of the expected
values, these terms contribute zero. Thus, the only contributions
to an odd moment come when each $b_j$ in the expansion occurs at
least twice, and at least one occurs three times. Hence, if $k =
2m+1$, we see we have at most $m+1$ degrees of freedom, this
coming from the case $b_{j_1}^3 b_{j_2}^2 \cdots b_{j_m}^2$. There
are $m$ different factors of $b$, and then we can choose any one
subscript. Once we have specified a subscript and which diagonals
we are on, the remaining subscripts are determined. As all moments
are finite, we find

\be M_{2m+1}(N) \ \ll_m \ \frac{1}{N^{\frac{2m+1}{2}+1}} N^{m+1} \
\ll_m \ \frac{1}{\sqrt{N}}. \ee

\subsection{Bounds for the Even Moments}\label{secboundseven}

We proceed in stages in calculating $M_{2k}(N)$, $2k \ge 4$.
First, we bound $M_{2k}(N)$ by $2^k \cdot 2^k \cdot (2k-1)!!$,
where $(2k-1)!!$ is the $2k^{th}$ moment of the Gaussian. We then
show that each factor of $2^k$ can be removed, and then show a
strict inequality holds.

\be M_{2k}(N) \ = \ \frac{1}{N^{k+1}}  \sum_{1 \leq i_{1},\cdots,
i_{2k} \leq N}\E(b_{|i_{1}-i_{2}|}b_{|i_{2}-i_{3}|}\dotsm
b_{|i_{2k}-i_{1}|}). \ee

If any $b_j$ occurs to the first power, its expected value is zero
and there is no contribution. Thus, the $b_j$s must be matched at
least in pairs. If any $b_j$ occurs to the third or higher power,
there are less than $k+1$ degrees of freedom, and there will be no
contribution in the limit.

The $b_j$s are matched in pairs, say $b_{|i_m - i_{m+1}|} =
b_{|i_n - i_{n+1}|}$. Let $x_m = |i_m - i_{m+1}| = |i_n -
i_{n+1}|$. There are two possibilities:

\be\label{eqchoicesigns} i_m - i_{m+1} \ = \ i_n - i_{n+1} \ \ \ \
\ \mbox{or} \ \ \ \ \ i_m - i_{m+1} \ = \ -(i_n - i_{n+1}). \ee

There are $k$ such pairs, thus we have $2^k$ choices of sign.
Further, there are $(2k-1)!!$ ways to pair off $2k$ numbers into
groups of two.

Fix a choice of sign and a pairing. Once we specify $x_1, \dots,
x_k$ and any one index, say $i_1$, all the other indices are
almost determined (if the choices are consistent). There is one
remaining freedom. After we've chosen which differences to match
and the values of these differences and the choice of signs, for
each time when there is a negative sign, there is one additional
choice: does the positive or negative difference occur first?
Thus, after we specify for each pair whether the positive or
negative difference occurs first, then all the indices are
determined.

Therefore, there are $N^{k+1}$ degrees of freedom. If all the
$x_j$s are distinct, we have the expected value of the second
moment of $p$, $k$ times. These contribute at most

\be \frac{1}{N^{k+1}} \cdot 2^k \cdot 2^k (2k-1)!! N^{k+1} \ = \
2^k (2k-1)!!. \ee

If some of the $x_j$s are equal, we have fewer than $k+1$ degrees
of freedom. We now have the expected value of a product of
moments of $p$, which is finite and independent of $N$. These
terms will not contribute in the limit. Therefore

\be \lim_{N\to \infty} \mtk \ \le \ 2^k \cdot 2^k (2k-1)!!. \ee

We now remove the factor of $2^k$ coming from the choice of signs.
Consider a pairing of the $b_j$s. We claim the only term which
contributes in the limit is when all signs are negative.

Let $x_1, \dots, x_k$ be the values of the $|i_j - i_{j+1}|$s, and
let $\gep_1, \dots, \gep_k$ be the choices of sign (see Equation
\ref{eqchoicesigns}). Define $\widetilde{x}_1 = i_1 - i_2$,
$\widetilde{x}_2 = i_2 - i_3, \dots, \widetilde{x}_{2k} = i_{2k} -
i_1$. Note exactly one $\widetilde{x}_j$ is $x_j$ and exactly one
is $\gep_j x_j$. We have

\bea i_2 & \ = \ & i_1 - \widetilde{x}_1 \nonumber\\ i_3 & = & i_1
- \widetilde{x}_1 - \widetilde{x}_2 \nonumber\\ & \vdots &
\nonumber\\ i_1 & = & i_1 - \widetilde{x}_1 - \cdots -
\widetilde{x}_{2k}. \eea

Therefore

\be \widetilde{x}_1 + \cdots + \widetilde{x}_{2k} \ = \
\sum_{j=1}^k (1 + \gep_j) x_j \ = \ 0. \ee

If any $\gep_j = 1$, then the $x_j$ are not linearly independent,
and we have fewer than $k+1$ degrees of freedom; these terms will
not contribute in the limit. Thus, the only valid assignment is to
have all the signs negative. There are now $2^k$ possible choices
of order (whether the negative or positive difference occurs
first), giving $2^k \cdot N^{k+1}$. We eliminate $2^k$ by changing
our viewpoint.

We have $k+1$ degrees of freedom. We match our differences into
$k$ pairs. Choose $i_1$ and $i_2$. We now look at the freedom to
choose the remaining indices $i_j$. Once $i_1$ and $i_2$ are
specified, we have $i_1-i_2$, and a later difference must be the
negative of that. If $i_2-i_3$ is matched with $i_1-i_2$, then
$i_3$ is uniquely determined (because it must give the opposite of
the earlier difference). If not, $i_3$ is a new variable. Now look
at $i_4$. If $i_3 - i_4$ is matched with an earlier difference,
then the sign of its difference is known, and $i_4$ is uniquely
determined; if this difference belongs to a new pair not
previously encountered, than $i_4$ is a new variable and free.
Proceeding in this way, we note that if we encounter $i_n$ such
that $i_{n-1}-i_n$ is paired with a previous difference, the sign
of its difference is specified, and $i_n$ is uniquely determined;
otherwise, if this is a difference of a new pair, $i_n$ is a free
variable, with at most $N$ choices. Thus we see there are at most
$N^{k+1}$ choices (note not all choices will work, as for example
the final difference $i_{2n} - i_1$ is determined before we get
there, because of earlier choices).

More explicitly, having $k+1$ degrees of freedom does not imply
each term contributes fully -- we will see there are Diophantine
obstructions which bound the moments away from the Gaussian's.
However, each pairing and choice of sign contributes at most
$N^{k+1}$, and we have shown

\be M_{2k}(N) \ \le \ (2k-1)!! + O_k\left(\fon\right). \ee

\subsection{The Fourth Moment}

The fourth moment calculation highlights the Diophantine
obstructions encountered, which bound the moments away from the
Gaussian.

\be\label{eqfourmom} M_{4}(N)=\Inv{N^{3}} \sum_{1 \leq
i_{1},i_{2},i_{3},i_{4}\leq N}
\E(\bs{1}{2}\bs{2}{3}\bs{3}{4}\bs{4}{1}) \ee

Let $x_{j}=|i_{j}-i_{j+1}|$. If any $b_{x_j}$ occurs to the first
power, its expected value is zero. Thus, either the $x_j$ are
matched in pairs (with different values), or all four are equal
(in which case they are still matched in pairs). There are 3
possible matchings; however, by symmetry (simply relabel), we see
the contribution from $x_1 = x_2$, $x_3 = x_4$ is the same as the
contribution from $x_1 = x_4$, $x_2 = x_3$.

If $x_1 = x_2$, $x_3 = x_4$, we have

\be i_1 - i_2 \ = \ -(i_2 - i_3) \ \ \ \mbox{and} \ \ \ i_3 - i_4
\ = \ -(i_4 - i_1). \ee

Thus, $i_1 = i_3$ and $i_2$ and $i_4$ are arbitrary. Using these
three variables as our independent degrees of freedom, we see
there are $N^3$ such quadruples. Almost all of these will have
$x_1 \neq x_3$, and contribute $\E(b_{x_1}^2 b_{x_3}^2) = 1$.
Given $i_1$ and $i_2$, $N-1$ choices of $i_4$ yield $x_1 \neq
x_3$, and one choice yields the two equal. Letting $p_4$ denote
the fourth moment of $p$, we see this case contributes

\be \frac{1}{N^3} \Big( N^2(N-1) \cdot 1 + N^2 \cdot p_4\Big) \ =
\ 1 - \frac{1}{N} + \frac{p_4}{N} \ = \ 1 + O\left(\fon\right).
\ee

The other possibility is for $x_1 = x_3$ and $x_2 = x_4$.
Non-adjacent pairing is what leads to Diophantine obstructions,
which decreases the contribution to the moment. Now we have

\be i_1 - i_2 \ = \ -(i_3 - i_4) \ \ \ \mbox{and} \ \ \ i_2 - i_3
\ = \ -(i_4 - i_1). \ee

This yields \be i_1 \ = \ i_2 + i_4 - i_3, \ \ i_1, i_2, i_3, i_4
\in \{1,\dots,N\}. \ee The fact that each $i_j \in \{1,\dots,N\}$
is what leads to the Diophantine obstructions. In the first case,
we saw we had three independent variables, and $N^3 + O(N^2)$
choices that were mutually consistent. Now, it is possible for
choices of $i_2, i_3$ and $i_4$ to lead to impossible values for
$i_1$. For example, if $i_2, i_4 \ge \frac{2N}{3}$ and $i_3 <
\frac{N}{3}$, we see $i_1 > N$. Thus, there are at most $(1 -
\frac{1}{27})N^3$ valid choices. This is enough to show the
Gaussian moment is strictly greater; later we will see that if
there is one moment less than the Gaussian, all larger even
moments are also smaller.

The following lemma shows this case contributes $\frac{2}{3}$ to
the fourth moment.

\begin{lem} Let $I_N = \{1,\dots, N\}$.  Then  $\#\{x,y,z \in I_N:
1  \leq  x+y-z \leq  N\} = \frac{2}{3}N^3 + \frac{1}{3}N$.
\end{lem}

\begin{proof}

Say $x+y = S \in \{2, \dots, 2 N \}$. For $2 \leq S \leq N$, there
are $S-1$ choices of $z$, and for $S \geq N+1$, there are
$2N-S+1$. Similarly, the number of $x,y \in I_N$ with $x+y=S$ is
$S-1$ if $S \le N+1$ and $2N-S+1$ otherwise. The number of triples
is

\bea \sum_{S=2}^N (S-1)^2 + \sum_{S=N+1}^{2N} (2N - S + 1)^2 \ = \
\frac{2}{3}N^3 + \frac{1}{3}N. \eea

\end{proof}

Collecting all the pieces, we have shown

\begin{thm}[Fourth Moment] Let $p_4$ be the fourth moment of $p$.
Then \be M_4(N) \ = \ 2\frac{2}{3} + \frac{2(p_4 - 1)}{N} +
\frac{1}{N^2}. \ee
\end{thm}

\subsection{Sixth and Eight Moments}

Any even moment can be explicitly determined by brute-force
calculation, though deriving exact formulas as $k \to \infty$
requires handling involved combinatorics. To calculate the higher
moments, consider $2k$ points on the unit circle, and look at how
many different shapes we get when we match in pairs. We find
$M_{6}(N)=11$ (compared to the Gaussian's $15$), and $M_{8}(N)=64
\frac{4}{15}$ (compared to the Gaussian's $105$). For the sixth
moment, there are five different configurations:

\begin{center}
\scalebox{.45}[.45]{\includegraphics{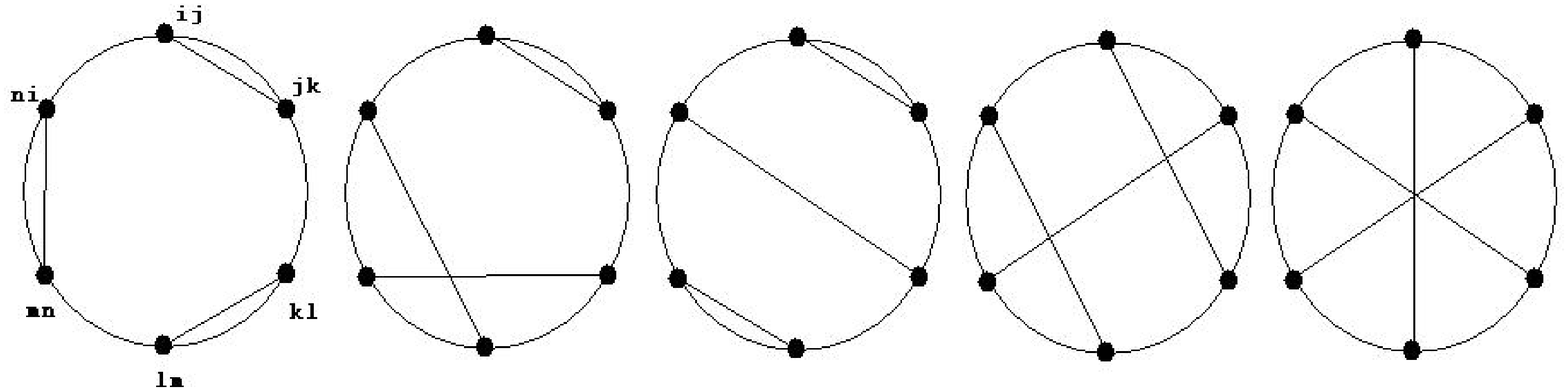}}
\end{center}

These occur $2, 6, 3, 3$ and $1$ time, contributing $1,
\frac{2}{3}, 1, \foh$, and $\foh$ (respectively); these correspond
to the $15 = (6 - 1)!!$ pairings. For the eight moment, the
smallest contribution is $\frac{1}{4}$, coming from the matching
$x_1 = x_3$, $x_2 = x_4$, $x_5 = x_7$, $x_6 = x_8$. It seems the
more crossings (in some sense), the greater the Diophantine
obstructions and the smaller the contribution.

\section{Upper Bounds of High Moments}

\subsection{Weak Upper Bound of High Moments}

\begin{lem} For $2k \ge 4$, $\lim_{N\to \infty} M_{2k}(N) <
(2k-1)!!$. \end{lem}
\begin{proof} Once we find a pairing that contributes less than
$1$ for some moment, we note that it will lift to pairings for
higher moments that will also contribute less than $1$. Say we
have such a pairing on $b_{|i_1-i_2|} \cdots b_{|i_{2k_0} - i_1|}$
giving less than 1. We extend this to a pairing on $2k > 2k_0$ as
follows. We now have

\be b_{|i_1-i_2|} \cdots b_{|i_{2k_0 - 1} - i_{2k_0}|}
b_{|i_{2k_0} - i_{2k_0+1}|} b_{|i_{2k_0+1} - i_{2k_0+2}|} \cdots
b_{|i_{2k-1} - i_{2k}|} b_{|i_{2k} - i_i|}. \ee

In groups of two, pair adjacent neighbors from $b_{|i_{2k_0+1} -
i_{2k_0+2}|}$ to $b_{|i_{2k-1} - i_{2k}|}$. This implies $i_{2k_0}
= i_{2k_0+2} = \cdots = i_{2k}$. Thus, looking at the first $2k_0
- 1$ and the last factor gives

\be b_{|i_1-i_2|} \cdots b_{|i_{2k_0 - 1} - i_{2k_0}|} b_{|i_{2k}
- i_i|} \ = \ b_{|i_1-i_2|} \cdots b_{|i_{2k_0 - 1} - i_{2k_0}|}
b_{|i_{2k_0} - i_i|}. \ee

Now pair these as in the pairing which gave less than $1$, and we
see this pairing will contribute less than $1$ as well.
\end{proof}

\subsection{Strong Upper Bound of High Moments}

In general, the further away one moment is from the Gaussian, the
more one can say about higher moments. While we do not have exact
asymptotics, one can show

\begin{thm} $\lim_{n\to \infty} \frac{M_{2k}}{(2k-1)!!} = 0$.
\end{thm}

\begin{proof} We will show that for any positive integer $c$, for $k$
sufficiently large, as $N \to \infty$ the moment is bounded by
$(\frac{2}{3})^c (2k-1)!!$. We have shown that we may take as
independent variables the $k$ values of the subscripts of the
$b_j$s ($x_1, \dots, x_k$) and any index. The goal is to show that
almost all of the pairings, for $k$ large, have at least $c$
Diophantine obstructions (of the type encountered in the fourth
moment). If there were no obstructions, these terms would
contribute $N^3$; the obstructions reduce the contribution to
$\frac{2}{3}N^3$.

We strategically replace our set of independent variables $i_d,
x_1, \dots, x_k$ with new variables which exhibit the
obstructions. We give full details on dealing with one
obstruction, and sketch how to add more. For simplicity, instead
of referring to $i_1, i_2, \dots, i_{2k}$, we use $i,j,k,\dots$
and $p,q,r,\dots$. Thus, in the trace expansion we have terms like
$a_{i_1i_2} = b_{|i_1-i_2|}$; we refer to this point by $i_1i_2$
or by $ij$.

\begin{center}
\scalebox{.45}[.45]{\includegraphics{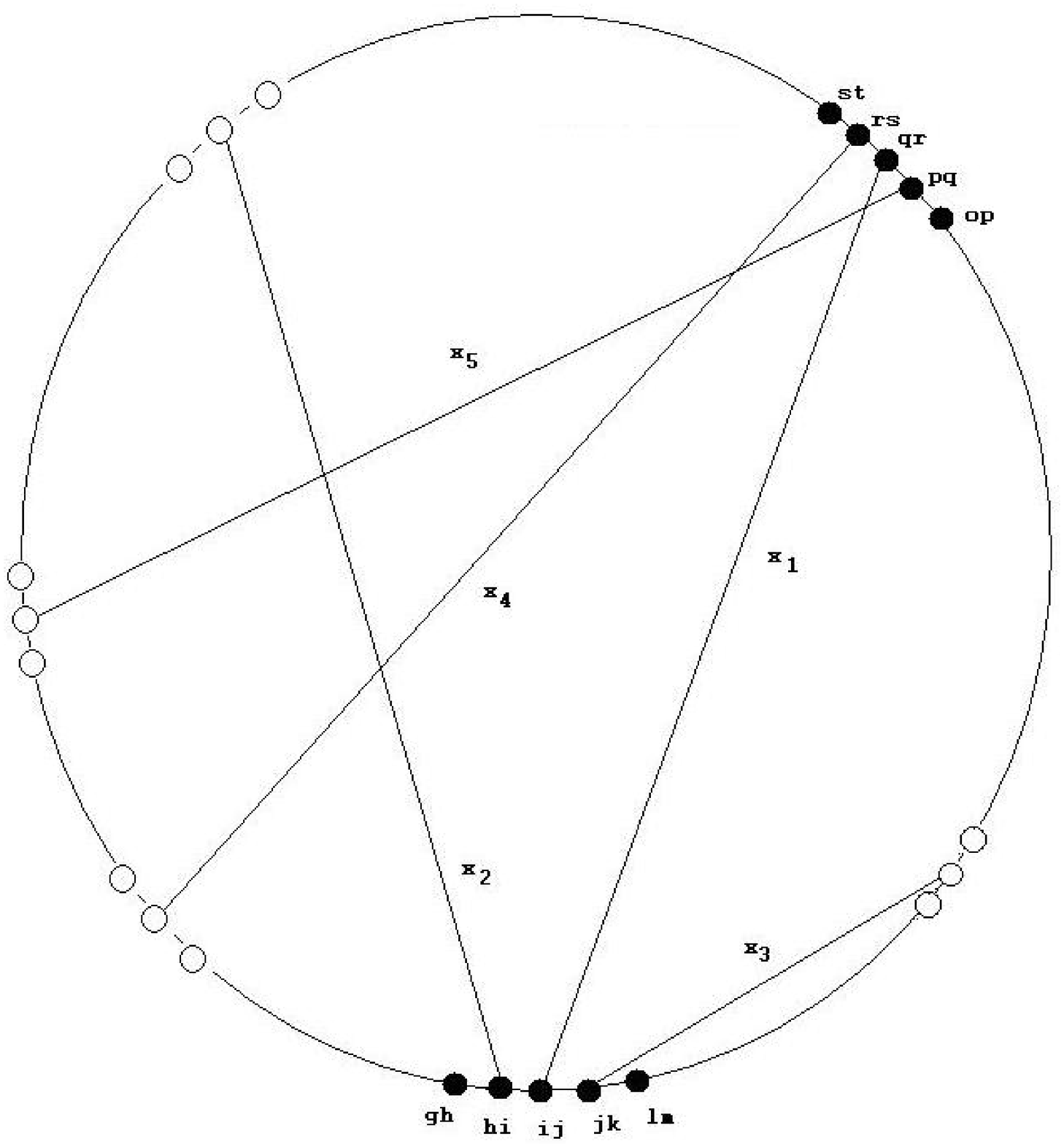}}
\end{center}

Say we pair $b_{|i-j|}$ with $b_{|q-r|}$. Let $x_1 = i - j =
-(q-r)$. If we knew $i = j + r - q$, with $j, r$ and $q$
independent free variables, then our earlier results show there
are only $\frac{2}{3}N^3$, not $N^3$, solutions. Unfortunately,
$j,r$ and $q$ need not be independent; however, for almost all of
the $(2k-1)!!$ pairings, they will be.

Create a buffer zone around $ij$ and $qr$ of two vertices on each
side, and assume that neither buffer zone intersects. Given $ij$,
there are $(2k-1) - 8$ possible choices to place $qr$. Now connect
the neighbors of $ij$ and $qr$ such that nothing is connected
within one vertex of another. There will be
$(2k-O(1))\cdot(2k-O(1))\cdot (2k-O(1))\cdot (2k-O(1))$ such
pairings. Note that, as we start placing some of these
connections, some vertices become unavailable. For example, say
there is exactly one vertex between the buffer of $ij$ and the
buffer of $qr$. This vertex is \emph{not} available for use, for
if we were to place another vertex there, the indices it gives
would not be independent. The same would be true if there were
just two vertices between the two buffers, and so on. In each
case, however, we only lose $O(1)$ vertices. As all these pairings
are separated, we may label their differences by $x_2, x_3, x_4$
and $x_5$, independent free variables.

The point is that the separation allows us to replace some the
independent variables $x_d$ with $j, r$ and $q$. Note that each
index appears in exactly two vertices on the circle, and they are
adjacent. Thus, these are the only occurrences of $i,j,q,r$ and we
may replace $x_5$ with $q$, $x_4$ with $r$, and $x_1$ with $j$. We
now have the desired situation: $i = j+r - q$, with all three on
the left independent free variables.

There are $(2k-11)!!$ ways to pair the remaining vertices. For
those pairs that have $j,q,r$ independent, the contribution is at
most $\frac{2}{3}N^3 \cdot N^{k+1-3}$; for the others, we bound
the contribution by $N^{k+1}$. Hence

\bea M_{2k}(N) & \ \le \ & \frac{1}{N^{k+1}} \left[(2k)^5
(2k-11)!! \frac{2}{3}N^{k+1} + O(k^4) \cdot (2k-11)!! \cdot
N^{k+1} \right] \nonumber\\ & \ \le \ & \frac{2}{3}(2k-1)!! +
O\left(\frac{(2k-1)!!}{k}\right). \eea

Therefore,

\be \frac{M_{2k}(N)}{(2k-1)!!} \ \le \ \frac{2}{3} + O\left(
\frac{1}{k} \right). \ee

There are two ways to handle the general case with $c$ Diophantine
obstructions. One may start with enormous buffer zones around the
initial pairs. As the construction progresses, we open up more and
more portions of the parts of the buffer zones not immediately
near the vertices. This keeps all but $O(1)$ vertices available
for use. Alternatively, along the lines of the first construction,
we can just note that by the end of stage $c$, $O_c(1)$ vertices
were unusable. We will still have the correct power of $2k$, with
a correction term smaller by a factor of $\frac{1}{k}$.
\end{proof}

\section{Lower Bound of High Moments}\label{seclowerboundmoments}

\subsection{Preliminaries}

By obtaining a sufficiently large lower bound for the even
moments, we show the limiting distribution has unbounded support.
In particular, we must find a lower bound $L_{2k}$ such that
$\lim_{k \to \infty} \sqrt[2k]{L_{2k}} = \infty$.

We know the moments are bounded by those of the Gaussian,
$(2k-1)!!$; the limiting value of the $2k$-th root of the Gaussian
(by Stirling's Formula) is $\frac{k}{e}$. We will show
$\sqrt[2k]{L_{2k}} \approx k^{\foh-\gep}$ in the limit.

The construction is as follows: in studying the $2k$-th moment, we
are led to sums of the form

\bea & & \frac{1}{N^{k+1}} \E\left[ \sum_{i_1=1}^N \cdots
\sum_{i_{2k}=1}^N a_{i_1,i_2} a_{i_2,i_3} \cdots
a_{i_{2k},i_1}\right] \nonumber\\ & \  = \ &\frac{1}{N^{k+1}}
\E\left[ \sum_{i_1=1}^N \cdots \sum_{i_{2k}=1}^N b_{|i_1-i_2|}
b_{|i_2-i_3|} \cdots b_{|i_{2k}-i_1|}\right]. \eea

If any $b_{|i_n-i_{n+1}|}$ occurs only once, as it is drawn from a
mean zero distribution, there is no contribution to the expected
value. Thus, the $2k$ numbers (the $b$s) are matched in at least
pairs, and, to obtain a lower bound, it is sufficient to consider
the case where the differences are matched in $k$ pairs. Let these
positive differences (of $|i_n - i_{n+1}|$) be $x_1, \dots, x_k$.

In Section \ref{secboundseven}, we showed the matchings must occur
with negative signs. Thus, if $|i_n - i_{n+1}| = |i_y-i_{y+1}|$,
then $(i_n-i_{n+1}) = -(i_y-i_{y+1})$. We let $\widetilde{x}_j =
i_j - i_{j+1}$. Thus, for any $x_j$, there is a unique $j_1$ such
that $\widetilde{x}_{j_1} = x_j$, and a unique $j_2$ such that
$\widetilde{x}_{j_2} = -x_j$. We call the first set of differences
positive, and the other set negative; we often denote these by
$\widetilde{x}_p$ and $\widetilde{x}_n$, and note that we have $k$
of each.

We have $k+1$ degrees of freedom. We may take these as the $k$
differences $x_k$, and then any index, say  $i_1$. We have the
relations

\bea\label{eqindicesms} i_2 & \ = \ & i_1 - \widetilde{x}_1 \nonumber\\
i_3 & = & i_1 - \widetilde{x}_1 - \widetilde{x}_2 \nonumber\\ &
\vdots & \nonumber\\ i_{2k} & = & i_1 - \widetilde{x}_1 - \cdots -
\widetilde{x}_{2k}. \eea

Once we specify $i_1$ and the differences $\widetilde{x}_1$
through $\widetilde{x}_{2k}$, all the indices are determined. If
everything is matched in pairs and each $i_j \in \{1,\dots,N\}$,
then we have a valid configuration, which will contribute $+1$ to
the $2k$-th moment. The reason it contributes $+1$ is because, as
everything is matched in pairs, we have the expected value of the
second moment of $p(x)$, $k$ times.

Thus, we need to show the number of valid configurations is
sufficiently large. The problem is that, in Equation
\ref{eqindicesms}, each index $i_j \in \{1,\dots, N\}$; however,
it is possible that a running sum $i_1 - \widetilde{x}_1 - \cdots
- \widetilde{x}_m$ is not in this range for some $m$. We will show
that we are often able to keep all these running sums in the
desired range.

\subsection{Construction}

Let $\alpha \in (\foh,1)$. Let $I_A = \{1,\dots, A\}$, where $A =
\frac{N}{k^\alpha}$. Choose each difference $x_j$ from $I_A$;
there are $A^k$ ways to do this. In the end, we want to study
$k$-tuples such that no value is chosen twice. Note such tuples
are lower order, namely there are at most $\ncr{k}{2} A^{k-1}$
such tuples. This is $O(N^{k-1})$. As $i_1$ takes on at most $N$
values (not all values will in general lead to valid
configurations), we see tuples with repeated values occur at most
$O(N^k)$ times; as we divide by $N^{k+1}$, these terms will not
contribute for fixed $k$ as $N \to \infty$. Thus, with probability
one (as $N \to \infty$), we may assume the $k$ values $x_j$ are
distinct.

Let us consider $k$ distinct positive numbers (the $x_j$s) drawn
from $I_A$, giving rise to $k$ positive differences
$\widetilde{x}_p$s and $k$ negative differences
$\widetilde{x}_n$s. Let us make half of the numbers
$\widetilde{x}_1, \dots, \widetilde{x}_k$ positive (arising from
the $\widetilde{x}_p$s), and half of these numbers negative
(arising from the $\widetilde{x}_n$s). Call this the first block
(of differences).

Then, in the differences $\widetilde{x}_{k+1},\dots,
\widetilde{x}_{2k}$ (the second block), we have the remaining
differences. Note every positive (negative) difference in
$\widetilde{x}_1,\dots, \widetilde{x}_k$ is paired with a negative
(positive) difference in $\widetilde{x}_{k+1},\dots,
\widetilde{x}_{2k}$. Note we have not specified the order of the
differences, just how many positive (negative) are in the first
block / second block.

Note two different $k$-tuples of differences $x_j$ \emph{cannot}
give rise to the same configuration (if we assume the differences
are distinct). This trivially follows from the fact that the
differences specify which diagonal of the Toeplitz matrix the
$a_{i_mi_{m+1}}$s are on; if we have different tuples, there is at
least one diagonal with an entry on one but not on the other.

Let us assume we have chosen the order of the differences in the
first block, $\widetilde{x}_1, \dots, \widetilde{x}_k$. We look at
a subset of possible ways to match these with differences in the
second block. In the second block, there are $\frac{k}{2}$
positive (negative) differences $\widetilde{x}_p$
($\widetilde{x}_n$). There are $(\frac{k}{2})!$ ways to choose the
relative order of the positive (negative) differences. Note we are
\emph{not} giving a complete ordering of the differences in the
second block. There are $k! > (\frac{k}{2})!^2$ ways to completely
order. We are merely specifying the relative order among the
positive (negative) elements, and not specifying how the positive
and negative differences are interspersed.

Thus, the number of matchings, each of which contribute $1$,
obtainable by this method is at most

\be\label{eqnumberprevalid} N \cdot (A^k - O(A^{k-1})) \cdot
(k/2)!^2, \ee

where $N$ is from the possible values for $i_1$, $A^k -
O(A^{k-1})$ is the number of $k$-tuples of distinct differences
$x_j \in I_A$, and $(k/2)!^2$ is the number of relative
arrangements of the positive and negative differences in the
second block (each of which is matched with an opposite difference
in the first block).

Not all of the above will yield a $+1$ contribution to the $2k$-th
moment. Remember, each index $i_m$ must be in $\{1,\dots,N\}$. We
now show that for a large number of the above configurations, we
do have all indices appropriately restricted. We call such a
configuration \emph{valid}.

\subsection{Number of Valid Configurations}

Most of the time, the sum of the positive differences
$\widetilde{x}_p$ in the first block will be close to the negative
of the sum of the negative differences $\widetilde{x}_n$ in the
first block.

Explicitly, we may regard the $\widetilde{x}_p$s
($\widetilde{x}_n$s) as independent random variables taken from
the uniform distribution on $I_A$ ($-I_A$) with mean approximately
$\foh A$ ($-\foh A$) and standard deviation approximately
$\frac{1}{2\sqrt{3}}A$. By the Central Limit Theorem, for $k$
large, the sum of the $\frac{k}{2}$ positive (negative)
$\widetilde{x}_p$s ($\widetilde{x}_n$s) in the first block
converges to a normal distribution with mean approximately
$\frac{kA}{4}$ ($-\frac{kA}{4}$) and standard deviation
approximately $\sqrt{\frac{k}{2}} \cdot\frac{A}{2\sqrt{3}}$.

Thus, for $N$ and $k$ sufficiently large, the probability that the
sum of the positive differences in the first block is in
$[\frac{kA}{4} - \frac{\sqrt{k} A}{2\sqrt{6}}, \frac{kA}{4} +
\frac{\sqrt{k}A}{2\sqrt{6}}]$ is at least $\foh$ (and a similar
statement for the negatives). Thus, of the $A^k$ tuples, at least
$\frac{1}{4} A^k$ will have the sums of the positive (negative)
differences lying in this interval (in the negative of this
interval). We call such choices \emph{good}.

Remember, in the arguments leading up to Equation
\ref{eqnumberprevalid}, we only specified two items. First, the
\emph{absolute values} of the $k$ differences (all distinct);
second, that half the positive differences are in the first block,
and the \emph{relative} orderings of the positive (negative)
differences in the second block is given.

Thus, we have freedom to choose how to intersperse the positives
and negatives in the first and second blocks. Consider a good
choice of $x_k$s. We place these differences in the first block of
length $k$ as follows. Choose the first positive difference from
our good list, and make the first difference positive. Keep
assigning (in order) the positive differences from our good list
until the running sum of the differences assigned to the first
block exceeds $A$. Then assign the negative differences from our
good list until the running sum of differences in the first block
is less than $-A$. We then assign positive differences again until
the running sum exceeds $A$, and so on. We assign half the
positive (negative) differences to the first block.

Throughout the process, the largest the running sum can be in
absolute value is $\max(2A, 2 \cdot \frac{\sqrt{k}A}{2\sqrt{6}})$.
This is because the $\frac{k}{2}$ positive (negative) differences
yield sums whose negatives are very close to each other, and each
added difference can change the running sum by at most $\pm A$.

We now assign the differences in the second block. We have already
chosen the positive and negative differences. There are
$(\frac{k}{2})!$ orderings of the positive (negative) differences.
We choose these relative orderings, and now choose how to
intersperse these. We put down the differences, again making sure
the running sum never exceeds in absolute value $\max(2A, 2 \cdot
\frac{\sqrt{k}A}{2\sqrt{6}})$.

Let $i_1 = 0$. From Equation \ref{eqindicesms}, we now see that
each index is at most $2\max(2A, 2 \cdot
\frac{\sqrt{k}A}{2\sqrt{6}})$. Therefore, each index is in
$\left[-\frac{2}{\sqrt{6}} \frac{N}{k^{\alpha-\foh}},
\frac{2}{\sqrt{6}} \frac{N}{k^{\alpha-\foh}}\right]$. Thus, if we
shift $i_1$ so that $i_1 \in \left[\frac{7}{8}
\frac{N}{k^{\alpha-\foh}}, \frac{N}{k^{\alpha-\foh}}\right]$, as
$\alpha
> \foh$ for $k$ large all indices will now be in $\{1,\dots,N\}$.
Thus, this is a valid assignment of indices.

We now count the number of valid assignments. We see this is at
least

\be \left(\frac{1}{8} \frac{N}{k^{\alpha-\foh}}\right) \cdot
\left(\frac{1}{4} A^k - \ncr{k}{2} A^{k-1}\right) \cdot (k/2)!^2.
\ee

To calculate the contribution to the $2k$-th moment from this
pairing, we divide by $N^{k+1}$. If any of the differences are the
same, there is a slight complication; however, as $N$ is large
relative to $k$, we may remove the small number of cases (at most
$\ncr{k}{2}A^k$) when we have repeat differences among the
$\widetilde{x}_p$s and $\widetilde{x}_n$s. By Stirling's Formula,
the main term is

\be \frac{1}{N^{k+1}} \frac{1}{32} \frac{N^{k+1}}{k^{k\alpha -
\foh}} \left(e^{\frac{k}{2} \log \frac{k}{2} - \frac{k}{2}} \sqrt{
2\pi (k/2)}\right)^2 \ = \ \frac{\pi k^{\frac{3}{2}}}{16
e^{(1+\log 2)k}} \cdot e^{(1-\alpha)k\log k}. \ee

Thus, the $2k$-th root looks like $\frac{e^{(1-\alpha)\log
k}}{e^{1+\log 2}} > O(k^{1-\alpha})$, proving the support is
unbounded.

\section{Weak Convergence}

We need to show that the variances tend to 0. Thus, we must show

\be \lim_{N \to \infty} \Big(\E[ M_m(A,N)^2 ] - \E[M_m(A,N)]^2
\Big) \ = \ 0. \ee

As $M_m(A,N) = \frac{1}{N^{\frac{m}{2}+1}}\mbox{Trace}(A^m)$, we
have

\bea \E[M_m(A,N)^2] & \ = \ & \frac{1}{N^{m+2}} \sum_{1 \le
i_1,\dots,i_m \le N} \sum_{1 \le j_1,\dots,j_m \le N} \E[
b_{|i_1-i_2|} \cdots b_{|i_m-i_1|} b_{|j_1-j_2|} \cdots
b_{|j_m-j_1|} ] \nonumber\\ \E[M_m(A,N)]^2 & = & \frac{1}{N^{m+2}}
\sum_{1 \le i_1,\dots,i_m \le N} \E[ b_{|i_1-i_2|} \cdots
b_{|i_m-i_1|}] \sum_{1 \le j_1,\dots,j_m \le N} E[b_{|j_1-j_2|}
\cdots b_{|j_m-j_1|} ]. \nonumber\\ \eea

There are two possibilities: if the absolute values of the
differences from the $i$s are completely disjoint from those of
the $j$s, then these contribute equally to $\E[ M_m(A,N)^2 ]$ and
$\E[M_m(A,N)]^2$. We are left with estimating the difference for
the crossover cases, when the value of an $i_\alpha - i_{\alpha+1}
= \pm (j_\beta - j_{\beta + 1})$.

We assume $m = 2k$; a similar proof works for odd $m$. Note
$N^{m+2} = N^{2k+2}$. The following two lemmas imply the variance
tends to 0. As our moments grow slower than the Gaussian, we
satisfy the conditions necessary to obtain almost surely weak
convergence.

\begin{lem} The contribution from crossovers in $\E[M_{2k}(A,N)]^2$
is $O_k(\frac{1}{N})$. \end{lem}

\begin{proof} For $\E[M_{2k}(A,N)]$, the expected value vanishes if anything is
unpaired. Thus, in $\E[M_{2k}(A,N)]^2$, in the $i$s and $j$s
everything is at least paired, and there is at least one common
value from a crossover. The maximum number of such possibilities
occurs when everything is paired on each side, and just one set of
pairs crosses over; for this crossover there are $2$ ways to
choose sign. In this case, there are $k+1$ degrees of freedom in
the $i$s, and $k+1 - 1$ degrees of freedom in the $j$s (we lost
one degree of freedom from the crossover). Thus, these terms give
$O(N^{2k+1})$. Considering now matchings on each side with triple
or higher pairings, more crossovers, and the two possible
assignments of sign to the crossovers, we find that $i$s and $j$s
with a crossover contribute $O_k(\frac{1}{N})$ to
$\E[M_m(A,N)]^2$.
\end{proof}

\begin{lem} The contribution from crossovers in $\E[M_m(A,N)^2]$
is $O_k(\frac{1}{N})$. \end{lem}

\begin{proof} If neither the $i$
differences nor the $j$ differences have anything unpaired (ie,
everything is either paired or higher), and there is at least one
crossover, it is easy to see these terms are $O_k(\frac{1}{N})$.
The difficulty occurs when we have unmatched singletons on either
side. Assume there are unmatched differences among the $i$s. We
only increase the number of degrees of freedom by replacing triple
pairings and higher among the $i$s with pairs and singletons (note
we may lose these degrees of freedom as these must be crossed and
matched with the $j$s, but we can always cross these over to the
$j$s with no net loss of degrees of freedom). Similarly, we can
remove triple and higher pairings among the $j$s.

Assume there are $s_i > 0$ singletons and $k - \frac{s_i}{2}$
pairs on the $i$ side, $s_j \ge 0$ singletons on the $j$ side, and
$C \ge \max(s_i,s_j)$ crossings. Note $s_j$ can equal 0, if we
send the singletons on the $i$ side to matched pairs among the
$j$s, but $C$ cannot be less than $s_i$ and $s_j$. Note $s_i, s_j$
are even.

On the $i$ side, there are $1 + (k - \frac{s_i}{2}) + (s_i - 1)$
degrees of freedom; the $1$ is from the freedom of assigning any
value to one index, then we have $k - \frac{s_i}{2}$ from pairs,
and then the last singleton's value is determined, so we have just
$s_i-1$ additional degrees of freedom from singletons.

Assume $s_j > 0$. On the $j$ side, there could have been $1 + (k -
\frac{s_j}{2}) + (s_j - 1)$ degrees of freedom, but we know we
have $C$ crossings. This loses at least $C-1$ degrees of freedom
(it's possible the last, forced $j$ difference already equalled an
$i$ difference). Thus, the number of degrees of freedom is

\be \Big[1 + \left(k - \frac{s_i}{2}\right) + (s_i - 1)\Big] +
\Big[1 + \left(k - \frac{s_j}{2}\right) + (s_j - 1) - (C - 1)\Big]
\ = \ 2k + 1 - \foh(2C - s_i - s_j).  \ee

If $s_j = 0$, then there are $1 + k - C$ degrees of freedom on the
$j$ side, and we get $2k+1 - (C - \frac{s_i}{2})$ degrees of
freedom.

Thus, there are at most $2k+1$ degrees of freedom. Doing the
combinatorics for choices of sign and number of triples and higher
shows these terms contribute $O_k(\frac{1}{N})$.\end{proof}

\begin{thm} The measures $\mu_{A,N}(x)$ weakly converge to a
universal measure of unbounded support, independent of $p$.
\end{thm}

\begin{proof} As $M_k$ is less than the Gaussian's moments, the
$M_k$s uniquely determine a probability measure, which by Section
\ref{seclowerboundmoments} has unbounded support. As $\E[M_k(A,N)]
\to M_k$ and the variances tend to zero, standard arguments give
weak convergence.
\end{proof}

\section{Almost Sure Convergence}

\subsection{Expansions}

For convenience in presentation, we assume $p(x)$ is even (ie, the
odd moments vanish); we remark on the trivial modifications to
handle the additional book-keeping from general $p(x)$. We will
show

\be \lim_{N\to \infty} \E\left[ |M_m(A,N) - \E[M_m(A,N)]|^4
\right] \ = \ O\left(\frac{1}{N^2}\right). \ee

The above (plus Chebychev and Borel-Cantelli) will yield almost
sure convergence. Expanding this out, it is sufficient to study

\bea & & \E[ M_m(A,N)^4] - 4 \E[ M_m(A,N)^3 ] \E[ M_m(A,N)] + 6
\E[ M_m(A,N)^2] \E[ M_m(A,N)]^2 \nonumber\\ & \ & - 3 \E[
M_m(A,N)] \E[M_m(A,N)]^3. \eea

For even moments, we may write the pieces as

\bea \E[ M_{2m}(A,N)^4 ]& \ = \ & \frac{1}{N^{4m+4}} \sum_i \sum_j
\sum_k \sum_l \E[ b_{is} b_{js} b_{ks} b_{ls}] \nonumber\\
\E[ M_{2m}(A,N)^3] \E[M_{2m}(A,N)] & = &  \frac{1}{N^{4m+4}}
\sum_i \sum_j \sum_k \sum_l \E[ b_{is} b_{js} b_{ks}]
\E[b_{ls}],\nonumber\\ \eea

(note we combined the $\left( {4 \atop 3}\right)$ and $\left( {4
\atop 4} \right)$ terms) and so on, where for instance

\be E_1 \ = \ \E[ b_{is} b_{js} b_{ks} b_{ls} ] \ = \ \E[ b_{|i_1
- i_2|} \cdots b_{|j_{2m} - j_1|} b_{|k_1 - k_2|} \cdots
b_{|k_{2m}-k_1|} b_{|l_1-l_2|}\cdots b_{|l_{2m}-l_1|}]. \ee

We fix some notation. Denote the expected value sums above by
$E_1, E_2, E_3$ and $E_4$ (which occur with factors of $1, -4, 6$
and $-3$ respectively). For $h \in \{i,j,k,l\}$, let $b_h$ refer
to the differences in $b_{|h_1 - h_2|}\cdots b_{|h_{2m}-h_1|}$ If
a difference in a $b_h$ is matched with another difference in
$b_h$, we say this is an \emph{internal} matching; otherwise, it
is an \emph{external} matching. By a singleton, pair, triple,
quadruple and so on, we refer to matchings within a $b_h$ (ie, an
internal matching). Thus, a triple occurs when exactly three of
the differences in a $b_h$ are equal.

Let $p_a$ denote the $a$-th moment of $p(x)$. Note $p_2 = 1$. In
$\sum \E[  b_i b_j b_k b_l]$, if we have all differences occurring
twice, except for two different differences occurring four times
(two quadruples) and another different one occurring six times
(one sextuple), we would have $1^{2m-7} p_4^2 p_6$.

Note there are at most $4m+4$ degrees of freedom -- everything
must be matched in at least pairs (we have $8m$ total differences,
as we are looking at the fourth power of the $2m$-th moment), and
then each $b_h$ has at most one more degree of freedom (can choose
any index). Thus, any terms with a loss of at least two degrees of
freedom contribute at most $O(\frac{1}{N^2})$.

\subsection{Only Pairs and Singletons}

We show there is no net contribution if there are no triples or
higher, and then deal with that case afterwards.

\begin{lem} Assume in addition there are no singletons. Then the contribution
is $O(\frac{1}{N^2})$. \end{lem}

\begin{proof} If there are no matchings between $b_h$s, then
everything is independent, and we get $1 - 4 + 6 - 3 = 0$. If two
pairs are matched, we lose one degree of freedom. There are
$\ncr{4}{2} = 6$ ways to choose two out of $i,j,k,l$ to share a
match.

For the four expected value sums, we get the following
contributions: $\ncr{4}{2} p_4$ from $E_1$; $\ncr{3}{2}p_4 +
(6-\ncr{3}{2})$ from $E_2$ (three times the two pairs are in the
expected value of a product together, giving $p_4$; the other
three times they are separated, giving $p_2 = 1$); $\ncr{2}{2}p_4
+ (6-\ncr{2}{2})$ from $E_3$ (only once are the matched pairs
together); $\ncr{4}{2}$ from $E_4$. Combining yields

\be 1 \cdot 6p_4 - 4(3p_4 + 3) + 6(p_4 + 5) - 3(6) \ = \ 0. \ee

If at least three pairs are matched together, or two sets of two
pairs are matched together, we lose at least $2$ degrees of
freedom, giving a contribution of size $O(\frac{1}{N^2})$.
\end{proof}

The following lemmas are the cornerstone of the later
combinatorics:

\begin{lem}\label{lemonesing} If there is a singleton in $b_h$ paired with something in $b_g$,
then there is a loss of at least one degree of freedom.
\end{lem}

Note if every difference in a $b_h$ (all singletons) is paired
with a difference in $b_g$ (all singletons), we have a loss of one
degree of freedom. We can choose any index and $2m-1$ differences
in $b_h$; the last difference is now determined. Once we choose
one index in $b_g$, all other indices are determined, for a total
of $1 + (2m-1) + 1$ (instead of $2m+2$) degrees of freedom. Thus,
instead of being able to choose $2m$ differences freely, we could
only choose $2m-1$.

Note the above argument holds if instead of all singletons, we
have elements of $b_g$ and $b_h$ only matched internally and
externally with each other.

\begin{proof} As we can cycle the labels, we may assume that
$b_{|h_{2m}-h_1|}$ is the singleton. Note that once any index and
the values of the other differences are given, then $|h_{2m}-h_1|$
is determined. We would like to conclude it is not free, and we
have lost a degree of freedom.

Its value is forced, and it must equal the difference from another
$b_g$ ($h \neq g \in \{i,j,k,l\}$), say $b_{|g_a - g_{a+1}|}$. If
$b_{|g_a - g_{a+1}|}$ wasn't forced, we have just lost a degree of
freedom; if it was forced, then we have already lost a degree of
freedom. \end{proof}

\begin{rek} In the above, we did not need the matching to be with a
singleton -- a pair, triple or higher would also have worked.
\end{rek}

\begin{lem}\label{lemthreesing} If at least three of the $b_h$s have a singleton,
there is a loss of at least two degrees of freedom. \end{lem}

\begin{proof} If there is a matching of singletons from say $b_i$
and $b_j$, and another matching from $b_k$ and $b_l$, the lemma is
clear from above. Without loss of generality, the remaining case
is when a singleton from $b_i$ is matched with one from $b_j$, and
another singleton from $b_i$ is matched with one from $b_l$. We
then apply the previous lemma to $(b_j,b_i)$ and $(b_k,b_i)$.
\end{proof}

We can now prove

\begin{thm} The contribution when there are no triple or higher
internal pairings is at most $O(\frac{1}{N^2})$. \end{thm}

\begin{proof} It is sufficient to show the non-zero contributions all
lost at least two degrees of freedom. We have already handled the
case when there are no singletons. If three or four $b_h$s have a
singleton, we are done by Lemma \ref{lemthreesing}. If exactly two
have singletons, then there is no contribution in the $E_1$
through $E_4$, except for the cases when they are under the
expected value together (remember the mean of $p$ vanishes).

We have already lost a degree of freedom in this case; if any pair
in any $b_h$ is matched with a pair in a $b_g$, we lose another
degree of freedom. Thus, we may assume there are no matches with
four or more elements. Thus, every difference that occurs, occurs
exactly twice.

There are $\ncr{4}{2}=6$ ways to choose which two of the four
$b_h$s have singletons paired. The contribution from $E_1$ is $6$,
from $E_2$ is $3$ (3 of the 6 times they are under the expected
value together; the other 3 times they are separated, and the
expected value of a difference occurring once is 0), from $E_3$ is
$1$ (only $1$ of the 6 ways have them under the expected value
together), and from $E_4$ is $0$. Thus, we have a contribution of

\be 1 \cdot 6 - 4 \cdot 3 + 6 \cdot 1 - 3 \cdot 0 \ = \ 0. \ee

We are left with the case when the only singletons are in one
$b_h$. As we are assuming there are no triple or higher internal
matchings, these singletons must then be matched with pairs,
giving external triples; as the odd moments of $p(x)$ vanish,
there is no net contribution. \end{proof}

\begin{rek} If we do not assume the odd moments of $p$ vanish,
additional book-keeping yields the contribution is of size
$\frac{1}{N^2}$. If exactly two of the $b_h$s have singletons,
then each has at least two; we've already handled the case when
they are matched together. As no difference can be left unmatched,
we just need to study the case when we get four triples or two
triples and a pair; each clearly loses two degrees of freedom;

We are left with the case when only one $b_h$ has singletons. We
are down one degree of freedom already, so there cannot be another
non-forced matching. If there are at least four singletons, we are
done. If there are two singletons, we get two triples (either with
the same or different $b_g$s). Similar arguments as before yield
the contributions are

\be 1 \cdot 6p_3^2 - 4 \cdot 3p_3^2 + 6\cdot p_3^2 - 3\cdot 0 \ =
\ 0 \ee

if the two external triples involve matchings from $b_h$ to the
same $b_g$, and

\be 1 \cdot 4p_3^2 - 4 \cdot 3p_3^2 + 6 \cdot 0 - 3 \cdot 0  \ = \
0. \ee

\end{rek}

\subsection{Eliminating Triple and Higher Matchings}

\begin{lem} If there are no crossovers, there is no net
contribution. \end{lem}

\begin{proof} If there are no crossovers, the expected value of
the products are the products of the expected values. Thus, each
term becomes $\E[M_{2m}(A,N)]^4$, and $1-4+6-3 = 0$. \end{proof}

\begin{lem} If there are at least two triples among all of the
$b_h$s, the contribution is $O\left(\frac{1}{N^2}\right)$.
\end{lem}

\begin{proof} Everything must be matched in at least pairs (or its expected
value vanishes). If there are only two values among six
differences, then instead of getting 3 degrees of freedom, we get
$2$. This is enough to see decay like $O\left(\frac{1}{N}\right)$.
If we didn't assume $p(x)$ were even, we would have more work; as
the odd moments vanish, however, the two triples must be paired
with other differences, or with each other. In either case, we
lose at least one degree of freedom from each, completing the
proof.
\end{proof}

\begin{rek} Similarly, one can show there cannot be a triple and
anything higher than a triple. Further, we cannot have two
quadruples or more, as a quadruple or more loses one degree of
freedom (a quadruple is two pairs that are equal -- instead of
having two degrees of freedom, we now have one). \end{rek}

\begin{lem} If there is a quadruple, quintuple, or higher matchings
within a $b_h$, the contribution is $O(\frac{1}{N^2})$. \end{lem}

\begin{proof} There can be no sextuple or higher, as this gives at
least three pairs matched, yielding one degree of freedom (instead
of three). If there is a quadruple or quintuple, everything else
must be pairs or singletons. As the odd moments vanish, a
quintuple must be matched with at least a singleton, again giving
six points matched, but only one degree of freedom.

We are left with one quadruple (which gives a loss of one degree
of freedom) and all else pairs and singletons. No pairs can be
matched to the quadruple or each other, as we would then lose at
least two degrees of freedom. If there are any singletons, by
Lemma \ref{lemonesing} there is a loss of a degree of freedom. If
we have a quintuple or higher, this is enough to lose two degrees
of freedom. Thus, we need only study the case of all pairs and one
quadruple, with no external matchings.

As everything is independent, we find a contribution of

\be 1 \cdot p_4 - 4 \cdot p_4 + 6 \cdot p_4 - 3 p_4 \ = \ 0, \ee

where $p_4$ is the fourth moment of $p$. \end{proof}

\begin{lem} If there is only one triple (say in $b_h$), the contribution is
$O(\frac{1}{N^2})$. \end{lem}

\begin{proof} As odd moments vanish, the triple must be paired
with a singleton from another $b_h$; further, there must be at
least one singleton in the same $b_h$ as the triple (as there are
an even number of terms). We thus lose a degree of freedom from
the triple matched with a singleton (four points, but one instead
of two matches), and we lose a degree of freedom from the
singleton in the same $b_h$ as the triple (Lemma
\ref{lemonesing}). Thus, we have lost two degrees of freedom.
\end{proof}

We have proved

\begin{thm}\label{thmevenmoms} The contribution from having a triple or higher
internal matching is $O(\frac{1}{N^2})$. \end{thm}

\begin{rek} Similar arguments work for general $p(x)$. \end{rek}

\subsection{Odd Moments}

As the odd moments of $p(x)$ vanish, handling

\be \lim_{N\to \infty} \E\left[ |M_{2m+1}(A,N) -
\E[M_{2m+1}(A,N)]|^4 \right] \ = \ O\left(\frac{1}{N^2}\right) \ee

is significantly easier.

\begin{thm}\label{thmoddmoms} We lose at least two degrees of freedom above,
implying the expected value is $O\left(\frac{1}{N^2}\right)$.
\end{thm}

\begin{proof} In each $b_h$, there is at least one odd internal matching
(or singleton); thus, only $E_1$ can be non-zero. If there are
four (or more) internal triples (or higher), we lose at least two
degrees of freedom.

If there are exactly three internal triples, either two are
matched together and one is matched with a singleton, or all three
are matched with singletons; in both cases we lose at least two
degrees.

If there are exactly two internal triples, there must be at least
two $b_h$s with singletons. If the triples are matched with
singletons, we lose two degrees; if the triples are matched
together we lose one degree from that, and one more degree from
the singletons (Lemma \ref{lemonesing}).

If there is exactly one triple, at least three $b_h$s have
singletons, and similar arguments yield a loss of at least two
degrees.

If there are no triples, then by Lemma \ref{lemthreesing} there is
a loss of at least two degrees. \end{proof}

Combining Theorems \ref{thmevenmoms} and \ref{thmoddmoms} yields

\begin{thm}\label{thmprestrong}
\be \lim_{N\to \infty} \E\left[ |M_m(A,N) - \E[M_{2m+1}(A,N)]|^4
\right] \ = \ O\left(\frac{1}{N^2}\right). \ee
\end{thm}

\begin{rek} Similar arguments work for general $p(x)$. \end{rek}

\subsection{Almost Sure Convergence}

We show that we have almost sure convergence. We first introduce
some notation, and then show how this follows from Theorem
\ref{thmprestrong}.

Fix $p(x)$ as before. Let $\Omega_N$ be the outcome space $(T_N,
\prod_{i=1}^{N-1} p(b_i)db_i)$, where $T_N$ is the space of all $N
\times N$ Real Symmetric Toeplitz matrices. Let $\Omega$ be the
outcome space $(T_{\N}, \prod p)$, where $T_{\N}$ is the set of
all $\N \times \N$ Real Symmetric Toeplitz matrices and $\prod p$
is the product measure built from having the entries iidrv from
$p(x)$. For each $N$, we have projection maps from $\Omega$ to
$\Omega_N$. Thus, if $A \in T_{\N}$ is a Real Symmetric Toeplitz
matrices, then $A_N$ is the restriction obtained by looking at the
upper left $N \times N$ block of $A$.

We slightly adjust some notation from before. Let $\mu_{A_N}(x)dx$
be the probability measure associated to the Toeplitz $N\times N$
matrix $A_N$. Then

\bea \mu_{A_N}(x)dx & \ = \ & \frac{1}{N} \sum_{i=1}^N
\delta\left(x - \frac{\gl_i(A_N)}{\sqrt{N}}\right) \nonumber\\
M_m(A_N) & = & \int_{\R} x^m \mu_{A_N}(x)dx \nonumber\\ M_m(N) & =
& \E[M_m(A_N)] \nonumber\\ M_m & = & \lim_{N \to \infty} M_m(N).
\eea

As $N \to \infty$, $M_m(N)$ converges to $M_m$, and the
convergence for each $m$ is at the rate of $\frac{1}{N}$. The
expectation above is with respect to the product measure on $T_N$
built from $p(x)$.

We want to show that, for all $m$, as $N \to \infty$,

\be M_m(A_N) \ \longrightarrow \ M_m \ \ \ \mbox{almost surely}.
\ee

By the triangle inequality,

\be |M_m(A_N) - M_m| \ \le \ |M_m(A_N) - M_m(N)| + |M_m(N) - M_m|.
\ee

As the second term tends to zero, it is sufficient to show the
first tends to zero for almost all $A$.

Chebychev's Inequality states that for any random variable $X$
with mean zero and finite $m$-th moment that

\be \mbox{Prob}(|X| \ge \gep) \ \le \ \frac{\E[X^m]}{\gep^m}. \ee

Note $\E[ M_m(A_N) - M_m(N)] = 0$, and by Theorem
\ref{thmprestrong}, $M_m(A_N) - M_m(N)$ has finite fourth moment.
In fact, Chebychev's Inequality and Theorem \ref{thmprestrong}
yield

\be \mbox{Prob}(|M_m(A_N) - M_m(N)| \ge \gep) \ \le \
\frac{\E[|M_m(A_N) - M_m(N)|^4]}{\gep^4} \ \le \ \frac{C_m}{N^2
\gep^4}. \ee

The proof is completed by applying the following:

\begin{lem}[Borel-Cantelli] Let $B_i$ be a sequence of
events with $\sum_i \mbox{Prob}(B_i) < \infty$. Let

\be B \ = \ \left\{\omega: \omega \in \bigcap_{j=1}^\infty
\bigcup_{k=j}^\infty B_i\right\}. \ee

Then the probability of $B$ is zero. \end{lem}

In other words, an $\omega$ is in $B$ if and only if that $\omega$
is in infinitely many $B_i$, and the probability of events
$\omega$ which occur infinitely often is zero.

Fix a large $k$ and let

\be B_N^{(k,m)} \ = \ \{A \in T_\N: |M_m(A_N) - M_m(N)| \ge
\frac{1}{k}. \} \ee

We have seen that $\mbox{Prob}(B_N^{(k,m)}) \le \frac{C_m
k^4}{N^2}$. Thus, for fixed $m$ and $k$, the conditions of the
Borel-Cantelli Lemma are met, and we deduce that the probability
of $A \in T_N$ that occur in infinitely many $B_N^{(k,m)}$ is
zero. We now let $k \to \infty$, and find for any fixed $m$, as $N
\to \infty$, $M_m(A_N) \to M_m$ with probability one. Let
$B_m^{i.o.}$ be the probability zero sets where we do not have
such convergence.

Let $B^{i.o.} = \bigcup_{m=1}^\infty B_m^{i.o.}$. As a countable
union of probability zero sets has probability zero, we see that
$\mbox{Prob}(B^{i.o.}) = 0$; however, this is precisely the set
where for some $m$, we do not have pointwise convergence.

Thus, except for a set of probability zero, we find $M_m(A_N) \to
M_m$ for all $m$.

\section{Poissonian Behavior?}

As there are only $N-1$ degrees of freedom for the Toeplitz
Ensemble, and not $O(N^2)$, it is reasonable to believe the
spacings between adjacent normalized eigenvalues may differ from
those of full Real Symmetric Matrices. For example, band matrices
of width 1 are just diagonal matrices, and there the spacing is
Poissonian ($e^{-x}$); full Real Symmetric Matrices are
conjectured to have their spacing given by the GOE distribution
(which is well approximated by $Axe^{-Bx^2}$).

For $d$-regular graphs, there are $\frac{dN}{2}$ degrees of
freedom. It has been numerically observed (see \cite{JMRR} among
others) that the spacings between adjacent eigenvalues look GOE.

We chose 1000 Toeplitz matrices ($1000 \times 1000$), with entries
iidrv from the standard normal. We looked at the spacings between
the middle 11 normalized eigenvalues for each matrix, giving us 10
spacings. A plot of the spacings between normalized eigenvalues
looks Poissonian.

\begin{center}
\scalebox{.45}[.45]{\includegraphics{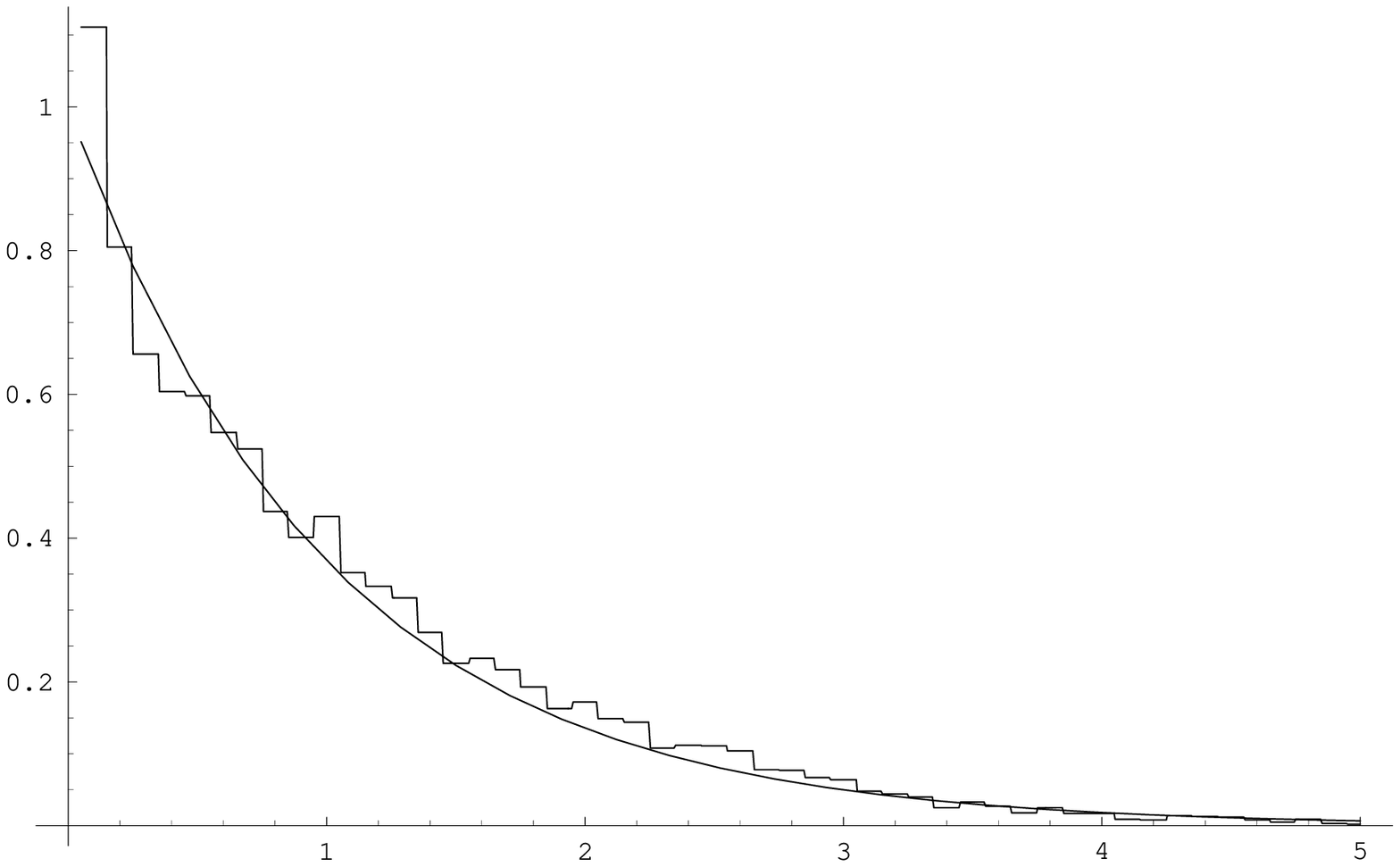}}
\end{center}

We conjecture that in the limit as $N \to \infty$, the local
spacings between adjacent normalized eigenvalues will be
Poissonian. It is interesting to note that Random $d$-Regular
Graphs have a comparable number of degrees of freedom; however, in
their adjacency matrices, there is significantly more independence
in the $a_{ij}$ -- for the Toeplitz Ensemble, we have a strict
structure, namely $a_{ij}$ depends only on $|i-j|$.

\end{document}